\documentclass[12pt]{article}

\newtheorem{theorem}{Theorem}

\newtheorem{proposition}[theorem]{Proposition}

\large

\title{{\bf  Regular spanning subgraphs of bipartite graphs of high minimum degree}}

\author{
B\'ela Csaba\thanks{Analysis and Stochastics Research Group of the Hungarian Academy of Sciences, 
University of Szeged, Hungary and Department of Mathematics, Western Kentucky University, Bowling Green, KY, USA}
\thanks{Partially supported by OTKA T049398 and by Hungarian State E\"otv\"os Fellowship.}
}

\date{}

\begin{document}
\maketitle

\begin{abstract}
Let $G$ be a simple balanced bipartite graph on $2n$ vertices, $\delta = \delta(G)/n$, and 
$\rho_0={\delta + \sqrt{2 \delta -1} \over 2}$. If $\delta \ge 1/2$ then $G$ has a
$\lfloor \rho_0 n \rfloor$-regular spanning subgraph. The statement is nearly tight. 
\end{abstract}

\section{Introduction}
In this paper we will consider regular spanning subgraphs of simple graphs. We mostly use standard graph 
theory notation: 
$V(G)$ and $E(G)$ will denote the vertex
and the edge set of a graph $G$, respectively. The degree of $x \in V(G)$ is denoted by 
$deg_G(x)$ 
(we may omit the subscript), $\delta(G)$ is the minimum degree of $G$.
We call a bipartite
graph $G(A,B)$ with color classes $A$ and $B$ {\it balanced} if $|A|=|B|$.
For $X, Y \subset V(G)$ we denote the number of edges of $G$ having one endpoint in $X$ and the other 
endpoint in $Y$ by $e(X,Y)$.
If $T \subset V(G)$ then
$G|_T$ denotes the subgraph we get after deleting every vertex of $V-T$ and the edges incident to them.
Finally, $K_{r, s}$ is the complete bipartite graph on color classes of size $r$ and $s$ for two positive 
integers $r$ and $s$. 

If $f: V(H) \rightarrow {\rm Z}^+$ is a function, then an $f$-factor is a subgraph $H'$ of the
graph $H$ such 
that $deg_{H'}(x) =f(x)$ for every $x \in V(H)$.
Notice, that when $f \equiv r$ for some $r \in {\rm Z}^+$, then $H'$ is an $r$-regular subgraph 
of $H$.

There are several results concerning $f$-factors of graphs. Perhaps the most notable among them is the
theorem of Tutte~\cite{T1}. Finding $f$-factors is in general not an easy task even for the case $f$ is
a constant and the graph is regular (see eg.,~\cite{BSW}).
In this paper we look for $f$-factors in (not necessarily regular) bipartite graphs with large 
minimum degree, for $f \equiv r$.

\begin{theorem}\label{tetel}
Let  $G(A,B)$ be a balanced bipartite graph on $2n$ vertices, and assume that
$\delta=\delta(G)/n \ge 1/2$. Set 
$\rho_0 =  {\delta+ \sqrt{2 \delta-1} \over 2}$. Then \\
(I) $G$ has a $\lfloor \rho_0 n \rfloor$--regular spanning subgraph;\\ 
(II) moreover, for every $\delta > 1/2$ if  
$n$ is sufficiently large and 
$\delta n$ is an integer then
there exists a balanced bipartite  graph  $G_{\delta}$ having minimum degree $\delta$ such that it does not admit a spanning regular subgraph
of degree larger than $\lceil \rho_0 n \rceil$.
\end{theorem} 

The above theorem plays a crucial role in the proof of some results 
in extremal graph theory (~\cite{CsM, Cs}).

\section{The main tool}

Let $F$ be a bipartite graph with color classes $A$ and $B$. By the well-known K\"onig--Hall theorem 
there is a perfect matching in $F$ if and only if
$|N(S)| \ge |S|$ for every $S \subset A$. We are going to need a far reaching generalization of this result, 
due to Gale and Ryser~\cite{GR,K} (one can find the proof in~\cite{L} as well). It gives a necessary and 
sufficient condition for the existence of an $f$-factor in a bipartite graph:

\begin{proposition} \label{faktor}
Let $F$ be a bipartite graph with bipartition $\{A,B \}$, and $f(x) \ge 0$ an integer valued function
on $A \cup B$. $F$ has an $f$--factor if and only if 
$$ (i) \sum_{x \in A} f(x) = \sum_{y \in B} f(y)$$
and \\
$$ (ii) \sum_{x \in X} f(x) \le e(X,Y) + \sum_{y \in B-Y} f(y)$$ \\
for all  $X \subset A$ and $Y \subset B$.
\end{proposition}

\section{Proof of Theorem~\ref{tetel}}

We will show the two parts of the theorem in separate subsections.

\subsection{Proof of part $I$}

Observe, that since we are looking for a spanning 
regular subgraph, the $f$ function of Proposition~\ref{faktor} will be identically $\rho n$ for
some constant $\rho$. 
We start with some notation: for $X \subset A$ let $\xi=|X|/n$, and for $Y \subset B$ let 
$\sigma=|Y|/n$.
We will normalize $e(X,Y)$: $\eta(X,Y)=e(X,Y)/n^2$. 
Let $$\eta_m(\xi,\sigma)=\min \{ \eta(X,Y): X \subset A, \ Y \subset B, \ |X|/n=\xi, |Y|/n=\sigma \}.$$ 

Since $f$ is identically $\rho n$, condition $(i)$ of Proposition~\ref{faktor} is satisfied. 
Moreover, if $\rho n$ is an integer and
$$\rho(\xi+\sigma-1) \le \eta_m(\xi, \sigma)$$
for some $\rho$  and for every $0 \le \xi, \sigma \le 1$, then $(ii)$ is satisfied, hence, $G$ has 
a $\rho n$--regular 
spanning subgraph. In the rest of this section we will show that the above inequality is valid for 
$\rho = \lfloor \rho_0 n \rfloor /n$.

Clearly, $e(X,Y) \ge |X|(\delta n -|B-Y|)$ and $e(X,Y) \ge |Y|(\delta n -|A-X|)$ for arbitrary sets 
$X \subset A$ and $Y \subset B$.  Hence, we have that $\eta_m(\xi,\sigma) \ge \max 
{ (\xi(\delta +\sigma -1), \sigma(\delta+ \xi -1)) }$. (In fact we always have that $\eta_m \ge 0$, since 
it is the edge density between the two color classes in $G$.)

\medskip

First consider the case $\xi = \sigma$. We are looking  for a $\rho$ for which 
$\rho (2 \xi -1)  \le \xi (\delta + \xi - 1)$. 
In another form, we need that 
$$p_\rho(\xi)=\xi^2 +(\delta-2\rho-1)\xi +\rho \ge 0.$$

The discriminant of the above polynomial is the polynomial $dcr(\rho)=
4 \rho^2- 4 \delta \rho +\delta^2-2 \delta +1$. Clearly, if $dcr(\rho) \le 0$ for some 
$\rho$, then $p_{\rho}(\xi) \ge 0$.

One can directly find the roots of $dcr(\rho)$:  ${\delta \pm \sqrt{2\delta-1} \over 2}$. 
At this point we have to be careful, since the degrees in a graph are non-negative
integers, so $\rho n$ has to be a natural number. 
We will show that $dcr(\rho) \le 0$ for $\rho=\lfloor (\delta + \sqrt{2 \delta -1})n/2 \rfloor /n$. 

\smallskip

Clearly, $dcr(x) \le 0$ in $I=[(\delta - \sqrt{2 \delta -1})/2, (\delta + \sqrt{2 \delta -1})/2]$,
the length of this interval is $\sqrt{2\delta -1}$. Divide the $[0,1]$ interval into $n$ disjoint subintervals
each of length $1/n$, denote the set of the endpoints of these subintervals by $S$. 
Observe that is $I \cap S \ge 1$, then we can pick the largest point of this intersection, this is $\rho \in I$, 
and we are done with proving that $p_{\rho}(\xi) \ge 0.$ 

We will investigate two cases: first, if $\delta > 1/2$, and second, if $\delta=1/2$.

\medskip

\noindent {\bf First case:} $\delta > 1/2$. We know that $\delta n$ is an integer, it is larger than $n/2$, hence, 
$\delta n \ge {n+1 \over 2}$. If the length of $I$ is at least $1/n$, it will intersect with $S$.
Assuming that $\sqrt{2 \delta -1} < {1 \over n}$ we would get $\delta n< {n^2+1 \over 2n}$,
but the latter expression is less than ${n+1 \over 2}$. Hence,
in this case $|I \cap S| \ge 1$. 

Let $g(\xi,\sigma)= \sigma(\delta+\xi-1)-\rho(\xi+\sigma-1)$. 
We will show, that $g(\xi, \sigma) \ge 0$ for $0 \le \sigma \le \xi \le 1$, in the lower right triangle $T$ 
of the unit square. This will prove that $(ii)$ of Proposition~\ref{faktor} is satisfied. 
Notice, that $g$ is bounded in the  triangle above, $-2 \le g(\xi, \sigma) \le 
\eta_m(\xi, \sigma)-\rho(\xi+\sigma-1)$, and continously differentiable.

Let us check the sign of $g$ on the border of the triangle. Since 
$\rho=\lfloor (\delta + \sqrt{2 \delta -1})n/2 \rfloor /n$, 
we have that $g(\xi, \xi) \ge 0$. $g(\xi,0)=-\rho(\xi-1) \ge 0$, and 
$g(1,\sigma)= \sigma (\delta - \rho) \ge 0$, because $\delta \ge (\delta + \sqrt{2 \delta -1})/2$.
Let us check the partial derivatives of $g$: 
$${\partial g \over \partial \xi}=\sigma-\rho,$$ and $${\partial g \over \partial \sigma}=\delta+\xi-1-\rho.$$
Assuming that $g$ achieves its minimum inside the triangle at the point $(\xi', \sigma')$ the partial 
derivatives of $g$ have to diminish at $(\xi',\sigma')$. It would then follow that 
$\sigma' = \rho$ and $\xi'=1+\rho-\delta$, therefore, $g(\xi',\sigma')=\rho^2 - \rho(2 \rho- \delta) = 
\delta \rho - \rho^2$. Hence $g$ is non-negative in $T$. 
The same reasoning works for the triangle $0 \le \xi \le \sigma \le 1$, this follows easily 
by symmetry. With this we finished the proof for the case $\delta > 1/2.$ 

\medskip

\noindent {\bf Second case:} $\delta = 1/2$. If $\delta n$ is even ($n$ is divisible by 4), we are done, 
since in this case $I$ contains the point
$\delta /2=\lfloor n/4 \rfloor /n$, and $\delta n/2$ is an integer. Therefore we have that $p_{1/4}(\xi) \ge 0$,
and as above, one can check that $g$ is non-negative in every point of $T$.

\smallskip

There is only one case left: if $\delta n$ is odd, that is, $n$ is of the form $4k+2$ for some natural
number $k$. In this case we want to prove, that that the spanning subgraph is ${k \over 4k+2}$-regular.

First observe, that for our purposes it is sufficient if $g(\xi,\sigma) \ge 0$ in a discrete point set: in the points 
$(\xi, \sigma)$ belonging to $(S \times S) \cap T$, since $|X|$ and $|Y|$ are natural numbers.  
Set $\rho={k \over 4k+2}$ and analyze the polynomial $p_{\rho}(\xi)$. It is an easy exercise to check that
it has two distinct roots: $1/2 ={2k+1 \over 4k+2}$ and $1/2 -1/n={2k \over 4k+2}.$ Hence, $p_{\rho}(\xi) \ge 0$
for $\xi \not\in (1/2-1/n, 1/2).$ 

We will cut out a small open triangle $T_s$ from $T$. $T_s$ has vertices
$(1/2,1/2),$ $(1/2, 1/2-1/n)$ and $(1/2-1/n, 1/2-1/n).$ Clearly, $T-T_s$ is closed and $T_s \cap (S \times S) =\emptyset$.

Recall, that $g(\xi, \sigma) = \sigma(\delta+\xi-1)-\rho(\xi+\sigma-1)$ for $0 \le \sigma \le \xi \le 1$. 
We will check the sign of $g$ on the 
border of $T-T_s$. There are two line segments for which we cannot apply our earlier results concerning $g.$
The first is 
$$L_1=\left\{(\xi,\sigma): {n-2 \over 2n} \le \xi \le {1 \over 2}, \sigma={n-2 \over 2n} \right\},$$
the second is
$$L_2=\left\{(\xi,\sigma): \xi={1 \over 2}, {n-2 \over 2n} \le \sigma \le {1 \over 2} \right\}.$$

On $L_1$ we get that 
$$g\left(\xi, {n-2 \over 2n}\right)= {n-2 \over 2n}\left(\xi - {1 \over 2}\right)-{k \over 4k+2}\left(\xi - {n-2 \over 2n}\right)=$$
$${k \over n}\xi - {k \over 2n}+ {k \over n^2}.$$
It is easy to see that the above expression is non-negative for every $(n-2)/(2n) \le \xi \le 1/2.$

For $L_2$ we have

$$g\left({1\over 2}, \sigma\right)=\sigma\left({1 \over 2}+{1 \over 2}-1\right) - {k \over 4k+2}\left({1 \over 2}+\sigma -1\right) =
{k \over 4k+2}\left({1 \over 2}-\sigma \right) \ge 0$$
for $(n-2)/(2n) \le \sigma \le 1/2.$

In order to finish proving that $g$ is non-negative in every point of $(S \times S) \cap (T-T_s)$ it is 
sufficient to show that the minimum of $g$ inside $T-T_s$ is at least as large as the minimum of $g$ 
on the border of $T-T_s$. This can be shown along the same lines as previously. By symmetry we will get 
that condition $(ii)$ of Proposition~\ref{faktor} is satisfied in every point of $S \times S.$

\subsection{Proof of part $II$}

For proving part $II$ of the theorem we want to construct a class of balanced
bipartite graphs the elements of which cannot have a large regular spanning subgraph.
We will achieve this goal in two steps. First, we will consider a simple linear function,
which, as we will see later, is closely related to our task. 
In the second step we will construct those bipartite graph which satisfy part $II$ of Theorem~\ref{tetel}.

Set $\gamma'={1-\sqrt{2\delta-1} \over 2}$ and let $0< p < 1.$
Consider the following equation:
 
$$(1-p)(1-\gamma') = \gamma'(1-p)+\delta -\gamma'. \quad \quad \quad \quad \quad \quad (1)$$

It is easy to see that $p'={\delta+\gamma'-1 \over 2\gamma' -1}$ is its solution. 
We have that

$$(1-p')(1-\gamma')=\left( 1-{\delta +\gamma' -1 \over 2\gamma' -1} \right)(1-\gamma')=
{\gamma' -\delta \over 2\gamma' -1}(1-\gamma').$$

Substituting $\gamma' ={1 -\sqrt{2 \delta -1} \over 2}$ we get 

$${\delta - {1 - \sqrt{2\delta-1} \over 2} \over \sqrt{2\delta-1}}
\left(1 - {1 -\sqrt{2\delta-1} \over 2}\right)=
{2\delta -1 + \sqrt{2\delta -1} \over \sqrt{2 \delta-1}}{1+\sqrt{2\delta -1} \over 2}=$$
$${1+\sqrt{2\delta -1} \over 2}{1+\sqrt{2\delta -1} \over 2}=
{\delta + \sqrt{2\delta-1} \over 2}.$$

\medskip

We promised to define a class of bipartite graphs for $\delta > 1/2$ which exist for every sufficiently 
large value of $n$ if $\delta n$ is a natural number, such that these graphs do not admit spanning 
regular graphs with large degree. 

For that let $\gamma=\lceil \gamma' n \rceil /n$. Then $\gamma n$ is an integer,
and $\gamma' \le \gamma \le \gamma' + 1/n.$ 
Let $G=(A,B,E)$ be a balanced bipartite graph on $2n$ vertices. $A$ is divided into two disjoint subsets, $A_l$ and
$A_e$, we also divide $B$ into $B_l$ and $B_e$. We will have that $|A_l|=|B_l|=\gamma n$ and $|A_e|=|B_e|=(1-\gamma)n$.
There are no edges in between the vertices of $A_l$ and $B_l$. 
The subgraphs $G|_{A_l \cup B_e}$ and  $G|_{B_l \cup A_e}$
are isomorphic to $K_{\gamma n, (1-\gamma)n}$, therefore, every vertex in $A_l \cup B_l$ has degree
$(1-\gamma)n$. We require that every vertex in $A_e \cup B_e$ has
degree $\delta n$, hence, $G|_{A_e \cup B_e}$ will be a $(\delta-\gamma)n$-regular graph.
Observe, that $\gamma < \delta < 1-\gamma$, thus, $\delta(G) = \delta n$.

Let us consider a simple method for edge removal from $G$: given $0< p <1$ discard 
$p(1-\gamma)n$ incident edges for every vertex in $A_l \cup B_l$, 
and no edge from $G|_{A_e \cup B_e}$. Of course, we need that $p(1-\gamma)n$ is an integer.

Then a vertex in $A_l \cup B_l$ will have degree $(1-p)(1-\gamma)n$, and the average 
degree of the vertices of $A_e \cup B_e$ will be $\gamma(1-p)n+(\delta - \gamma)n$.
Choose $\widetilde{p}$ to be the solution of the following equation:

$$(1-p)(1-\gamma)n =\gamma (1-p)n+(\delta - \gamma)n. \quad \quad \quad \quad \quad \quad (2)$$

Notice, that the only difference between (1) and (2) is that we subsituted $\gamma'$ by $\gamma$.
One can see that if $p < \widetilde{p}$ then there is a vertex $x \in A_e \cup B_e$ such that every vertex of 
$A_l \cup B_l$ will have degree larger than $deg(x)$. That is, for finding a regular subgraph more edges 
have to be discarded
among those which are incident to the vertices of $A_l \cup B_l$.

The solution of (2) is $\widetilde{p} ={\delta + \gamma -1 \over 2\gamma -1}$ (here $\widetilde{p} (1-\gamma)n$ is
not necessarily an integer). 
Computing the derivative shows 
that $\gamma \ge \gamma'$ implies $p' \ge \widetilde{p}$. Let us show that $p'-\widetilde{p}$ is small:

$$p'-\widetilde{p}={\delta+\gamma' - 1 \over 2\gamma' -1}-{\delta+\gamma -1 \over 2\gamma -1}=$$
$${(\delta+\gamma'-1)(2\gamma -1)-(\delta+\gamma-1)(2\gamma'-1) \over (2\gamma-1)(2\gamma'-1)}=
{2\gamma \delta - 2\gamma' \delta + \gamma' - \gamma \over (2\gamma -1)(2\gamma' -1)}.$$

Observe, that $1-2\gamma' =\sqrt{2\delta -1}$, and that $1-2\gamma \ge 1-2\gamma'- 2/n >0$ whenever $n$ is
sufficiently large. Therefore,

$${2\gamma \delta - 2\gamma' \delta + \gamma' - \gamma \over (2\gamma -1)(2\gamma' -1)}=
{(\gamma' - \gamma) \sqrt{2 \delta -1} \over 2\gamma -1} =$$
$$(\gamma - \gamma')
{\sqrt{2\delta -1} \over 1-2\gamma} \le {(\gamma - \gamma')\sqrt{2\delta -1} \over 1 -2\gamma' -2/n} \le $$
$${1 \over n} \left(1+{2 \over n\sqrt{2\delta -1}-2}\right)={1 \over n}\left(1+O({1/n})\right).$$

Above we used the fact that $\gamma \le \gamma' + {1 \over n}.$
Since $\widetilde{p} (1-\gamma)n$ is not necessarily an integer, we introduce $p_0$: 
$p_0=\lceil \widetilde{p}(1-\gamma)n \rceil /((1-\gamma)n).$ Clearly, the least number of edges one has to
remove from the vertices of $A_l \cup B_l$ in order to find a spanning regular subgraph of $G$ is at least
$p_0(1-\gamma)n$. With this choice of $p_0$ every degree in $A_l \cup B_l$ will be $(1-p_0)(1-\gamma)n$ after
the edge removal process. 

Finally, we show that $(1-p_0)(1-\gamma)$ is very close to ${\delta+\sqrt{2\delta-1} \over 2}:$

$$(1-p_0)(1-\gamma)-{\delta+\sqrt{2\delta-1} \over 2}=(1-p_0)(1-\gamma)-(1-p')(1-\gamma') \le $$
$$(1-\widetilde{p})(1-\gamma')-(1-p')(1-\gamma')=(1-\gamma')(1-\widetilde{p}-1+p')=$$
$$(1-\gamma')(p'-\widetilde{p})=(1-\gamma'){1 \over n}\left(1+O(1/n)\right).$$

If $n$ is sufficiently large, then $(1-\gamma')(1+O(1/n)) < 1$, since $0 < \gamma' < 1/2.$
Hence, if $H \subset G$ is an $r$-regular spanning subgraph, then 

$$\rho_0 n = \left\lfloor {\delta +\sqrt{2\delta -1} \over 2}n \right\rfloor \le r < (1-p')(1-\gamma')n +1 =  
{\delta +\sqrt{2\delta -1} \over 2}n +1.$$

Since $r$ is an integer which is less than ${\delta +\sqrt{2\delta -1} \over 2}n +1$,
we get that 

$$\left\lfloor {\delta +\sqrt{2\delta -1} \over 2}n \right\rfloor \le r \le 
\left\lceil {\delta +\sqrt{2\delta -1} \over 2}n \right\rceil,$$

and this is what we wanted to prove.


\bigskip




\begin{thebibliography}{20}

\bibitem{BSW} B. Bollob\'as, A. Saito, N. Wormald, {\it Regular factors of regular graphs}, Journal of Graph Theory 
{\bf 9} (1985) 97--103.

\bibitem{CsM} B. Csaba, M. Mydlarz, {\it Approximate multipartite version of the Hajnal--Szemer\'edi theorem},
submitted for publication.

\bibitem{Cs} B. Csaba, {\it Approximating the Bollob\'as--Eldridge--Catlin conjecture for bounded 
degree graphs}, manuscript.

\bibitem{K} M. Krause, {\it  A simple proof of the Gale-Ryser theorem}, The American Mathematical Monthly {\bf 4} (1996) 335--337.

\bibitem{L} L. Lov\'asz, {\it Combinatorial Problems and Exercises}, 
Akad\'emiai Kiad\'o, Budapest, 1993.

\bibitem{GR} H. J. Ryser, {\it Combinatorial mathematics}, Wiley and Sons, 1963

\bibitem{T1} W. T. Tutte, {\it The factors of graphs}, Canadian J. Math. {\bf 4} (1952) 314--328.


\end{thebibliography}
\end{document}